\documentclass[11pt,twoside]{article}
\usepackage{latexsym,amsmath}
\usepackage{indentfirst}
\usepackage{graphicx}

\numberwithin{equation}{section}
\newtheorem{Prop}{\bf Proposition}[section]

\newtheorem{Ex}{\bf Example}[section]

\begin{document}
\def \b{\Box}

\begin{center}
{\LARGE {\bf Stability analysis and control chaos for\\[0.1cm]
fractional $5D$ Maxwell-Bloch model}}\\
\end{center}

\begin{center}
{\bf Mihai IVAN}
\end{center}

\setcounter{page}{1}

\pagestyle{myheadings}

 {\small {\bf Abstract}. In this paper we investigate the dynamical behavior of fractional differential systems associated to
 $5D$ Maxwell-Bloch model in terms of fractional Caputo derivatives.}
{}

\section {Introduction}

The fractional calculus has been found to be an important tool in
various fields, such as mathematics, physics, engineering,
chemistry, biology, economics, chaotic dynamics, optimal control
and other complex dynamical systems \cite{ahme, podl, chen, giva,
mato}.

In this paper is used the Caputo definition of fractional
derivatives. Let $ f\in C^{\infty}( \textbf{R}) $ and $ \alpha \in
\textbf{R}, \alpha
> 0. $ The $\alpha-$order Caputo differential operator
\cite{difo}, is
described by\\[-0.2cm]
\[
D_{t}^{\alpha}f(t) = J^{m-\alpha}f^{(m)}(t), ~\alpha > 0,
\]
where $~f^{(m)}(t)$  represents the $ m-$order derivative of the
function $ f,~m \in \textbf {N}^{\ast}$ is an integer such that $
m-1 \leq \alpha \leq m $ and $ J^{\beta} $ is the $ \beta-$order
Riemann - Liouville integral operator \cite{podl}, which is
expressed by $~ J^{\beta}f(t)
=\displaystyle\frac{1}{\Gamma(\beta)}\int_{0}^{t}{(t-s)^{\beta
-1}}f(s)ds,~\beta > 0,$ where $~\Gamma $ is the Euler Gamma
function. If $ \alpha=1$, then $ D_{t}^{\alpha}f(t) =
\frac{df}{dt}.$

In this paper we suppose that $ \alpha \in (0,1].$

 The paper is structured as follows. In Section 2 we recall some results
concerning the study of stability for fractional systems. The
problem of the existence and uniqueness of solution for the
fractional $5D$ Maxwell-Bloch system $(3.3)$ is analyzed in
Section 3. Section 4 is devoted to studying of the stability of
equilibrium states for fractional system $(3.3)$. Also, the
unstable equilibrium states of this system can be controlled via
fractional stability theory. In Section 5, the numerical
integration and numerical simulation for the controlled fractional
$5D$ Maxwell-Bloch model $(4.1)$ are given.

\section{Preliminaries on fractional dynamical systems}

 We consider the following system of fractional differential equations on ${\bf R}^{n}$:\\[-0.2cm]
\begin{equation}
D_{t}^{\alpha}x^{i}(t) = f_{i}(x^{1}(t), x^{2}(t), \ldots,
x^{n}(t)) ,~~ i=\overline{1,n},\label{(2.1)}
\end{equation}
where $\alpha\in (0,1), f_{i}\in C^{\infty}({\bf R}^{n}, {\bf R}),
~ D_{t}^{\alpha} x^{i}(t)$ is the Caputo fractional derivative of
order $\alpha$ for $i=\overline{1,n}$ and $t\in [0,\tau)$ is the
time.

The fractional dynamical system $(2.1)$ can be written as follows:
\begin{equation}
D_{t}^{\alpha}x(t) = f(x(t)) ,\label{(2.2)}
\end{equation}
where $~f(x(t)) = (f_{1}(x^{1}(t),\ldots, x^{n}(t)),
f_{2}(x^{1}(t),\ldots, x^{n}(t)), \ldots, f_{n}(x^{1}(t),\ldots,
x^{n}(t)))^{T} $ and $ D_{t}^{\alpha} x(t)= ( D_{t}^{\alpha}
x^{1}(t), \ldots, D_{t}^{\alpha} x^{n}(t))^{T}.$

A point $x_{e}=(x_{e}^{1}, x_{e}^{2},\ldots, x_{e}^{n})\in {\bf
R}^{n}$ is said to be {\it equilibrium state} of the system
$(2.1)$, if $~D_{t}^{\alpha}x^{i}(t) =0$ for $i=\overline{1,n}$.

\markboth{M. Ivan}{ Stability analysis and control chaos for
fractional $5D$ Maxwell-Bloch model}

The equilibrium states of the fractional dynamical system $(2.1)$
are determined by solving the set of equations: $~f_{i}(x^{1}(t),
x^{2}(t), \ldots, x^{n}(t)) = 0 ,~~ i=\overline{1,n}.$

The Jacobian matrix associated to system $(2.1)$ is $~
J(x)=(\displaystyle\frac{\partial f_{i}}{\partial
x^{j}}),~~~i,j=\overline{1,n}.$

The stability of the system $(2.1)$ has been studied by Matignon
in \cite{mati}, where necessary and sufficient conditions have
been established.

\begin{Prop} {\rm (\cite{mati})}
Let $x_{e}$ be an equilibrium state of system $(2.1)$ and
$J(x_{e})$ be the Jacobian matrix $J(x)$ evaluated at $x_{e}$.

 $(i)~ x_{e}$ is
locally asymptotically stable, iff all eigenvalues
$\lambda(J(x_{e}))$ of  $ J(x_{e})$ satisfy:
\begin{equation}
| arg(\lambda (J(x_{e}))) | > \displaystyle\frac{\alpha
\pi}{2}.\label{(2.3)}
\end{equation}

 $(ii)~ x_{e}$ is
locally stable, iff either it is asymptotically stable, or the
critical eigenvalues of $ J(x_{e})$ which satisfy $~| arg(\lambda
(J(x_{e}))) | = \displaystyle\frac{\alpha \pi}{2}~$ have geometric
multiplicity one.\hfill$\Box$
\end{Prop}

In the case when $x_{e}$ is a unstable equilibrium state of the
fractional system $(2.2)$, we associate to $(2.2)$ a new
fractional system as follows.

 The {\it controlled fractional system} associated  to
system $(2.2)$ is described by:\\[-0.2cm]
\begin{equation}
D_{t}^{\alpha}x(t) = f(x(t))-k (x(t)-x_{e}) ,\label{(2.4)}
\end{equation}
where $~k= diag(k_{1},\ldots, k_{n}),~ k_{i}\geq 0, i =
\overline{1,n}$ and $x_{e}$ is an equilibrium state of $(2.2)$.

If one selects the appropriate parameters $k_{i},
i=\overline{1,n}$  which then make the eigenvalues of the
linearized equation of $(2.4)$ satisfy one of the conditions from
Proposition 2.1, then the trajectories of $(2.4)$ asymptotically
approaches the unstable equilibrium state $x_{e}$ in the sense
that $\lim_{t\rightarrow \infty} \|x(t)-x_{e}\|= 0$, where
$\|\cdot\|$ is the Euclidean norm.

\section{The fractional $5D$ Maxwell-Bloch model}

In the physics of self-induced transparency  for the most lasers
and the most atoms the so called two level lossless model is an
excellent approximation and is quite adequate for an understanding
of the basic physics behind many coherent transient phenomena
\cite{aleb}. Self-induced transparency equations based upon this
model are derived from the Maxwell-Schr$\ddot{o}$dinger equations
in the paper of Holm and Kovacic \cite{hoko}. More precisely,
after averaging and neglecting non-resonant terms, the unperturbed
Maxwell-Bloch dynamics in the rotating wave approximation (RWA)
can be written on ${\bf C}^{2}\times {\bf R}$ in the following
form:\\[-0.2cm]
\begin{equation}
\frac{du}{dt} = v,~~~ \frac{dv}{dt} = u w,~~~\frac{dw}{dt}
=\frac{1}{2}(\bar{u} v + u\bar{v}),\label{(3.1)}
\end{equation}
where the  superscript $"-"$ denotes the complex conjugation.
Physically speaking the complex scalar functions $u, v$ represent
the self-consistent electric field and respectively the
polarization of the laser-matter, the real scalar function $w$
describes the difference of its occupation numbers \cite{foho,
huan}.

Using the transformations $ u = x^{1}+ i x^{2}, ~ v = x^{3} + i
x^{4},~ w = x^{5},$ the dynamical system $(3.1)$ becomes:\\[-0.2cm]
\begin{equation}
 \dot{x}^{1} = x^{3},~~~ \dot{x}^{2} = x^{4},~~~
\dot{x}^{3} =  x^{1} x^{5},~~~
 \dot{x}^{4} = x^{2} x^{5},~~~
\dot{x}^{5} = -( x^{1} x^{3}+ x^{2} x^{4}),\label{(3.2)}\\[-0.1cm]
\end{equation}
where $ \dot{x}^{i}= \displaystyle\frac{x^{i}(t)}{dt}$  for
$i=\overline{1,5}.$ The phase space of  $(3.2)$ is ${\bf R}^{5}.$

 The dynamical system $(3.2)$ is called the {\it
five-dimensional Maxwell-Bloch equations} or the {\it $5D$
Maxwell-Bloch model}.

In \cite{foho}, Fordy and Holm discuss the phase space geometry of
the solutions of the system $(3.1)$ and show that it has three
Hamiltonian structures. More recently, Birtea and Ca\c su
\cite{bica} solve the stability problem for the isolated
equilibria of the system $(3.2)$.

The {\it fractional $5D$ Maxwell-Bloch model} associated to $5D$
Maxwell-Bloch model $(3.2)$ is defined by the following set of
equations:\\[-0.2cm]
\begin{equation}
\left\{ \begin{array}{ccl}
 D_{t}^{\alpha}{x}^{1}  &=&  x^{3}\\
  D_{t}^{\alpha}{x}^{2}  & =& x^{4}\\
 D_{t}^{\alpha}{x}^{3}  &=&  x^{1} x^{5},~~~~~~~~~~~~~~~~~\alpha\in (0,1).\\
  D_{t}^{\alpha}{x}^{4}  & =& x^{2} x^{5}\\
 D_{t}^{\alpha}{x}^{5}  &=& -( x^{1} x^{3}+ x^{2} x^{4})
\end{array}\right.\label{(3.3)}
\end{equation}

 The initial value problem of fractional model $(3.3)$ can be represented in the following
matrix
form:\\[-0.2cm]
\begin{equation}
D_{t}^{\alpha}x(t)  =  Ax(t) + x^{1}(t) A_{1} x(t) + x^{2}(t)
A_{2} x(t),~~~~x(0) = x_{0},\label{(3.4)}
\end{equation}
where $0 <\alpha < 1,~ x(t)= ( x^{1}(t),
 x^{2}(t), x^{3}(t), x^{4}(t), x^{5}(t) )^{T}, ~t\in(0,\tau)$ and\\[-0.1cm]
\[
A = \left ( \begin{array}{ccccc}
0 & 0 & 1 & 0 & 0\\
0 & 0 & 0 & 1 & 0\\
0 & 0 & 0 & 0 &  0\\
0 & 0 & 0 & 0 & 0\\
0 & 0 & 0 & 0 & 0\\
\end{array}\right ),~~ A_{1} = \left ( \begin{array}{ccccc}
0 & 0 & 0 & 0 & 0\\
0 & 0 & 0 & 0 & 0\\
0 & 0 & 0 & 0 & 1\\
0 & 0 & 0 & 0 & 0\\
0 & 0 &-1 & 0 & 0\\
\end{array}\right ),~~
A_{2} = \left ( \begin{array}{ccccc}
0 & 0 & 0 & 0 & 0\\
0 & 0 & 0 & 0 & 0\\
0 & 0 & 0 & 0 & 0\\
0 & 0 & 0 & 0 & 1\\
0 & 0 & 0 & -1 &0\\
\end{array}\right ).
\]
\begin{Prop}
The initial value problem of the fractional $5D$ Maxwell-Bloch
model $(3.4)$ has a unique solution.
\end{Prop}
{\it Proof.} Let $ f(x(t))= A x(t) + x^{1}(t) A_{1} x(t) +
x^{2}(t) A_{2} x(t).$ It is obviously continuous and bounded on $
D =\{ x \in {\bf R}^{5} |~ x^{1}\in [x_{0}^{1} - \delta, x_{0}^{1}
+ \delta],~ x^{2}\in [x_{0}^{2} - \delta, x_{0}^{2} + \delta]\} $
for any $\delta>0.$ We have $~f(x(t)) - f(x_{1}(t)) = A ( x(t) -
x_{1}(t)) + y(t) + z(t),~$ where $~y(t)=x^{1}(t) A_{1} x(t)-
x_{1}^{1}(t) A_{1} x_{1}(t) $
and $~z(t) = x^{2}(t) A_{2} x(t) - x_{1}^{2}(t) A_{2} x_{1}(t).$ Then\\[0.2cm]
$(1)~~|f(x(t)) - f(x_{1}(t))|\leq \|A\|\cdot |x(t)-x_{1}(t)| +
|y(t)| + | z(t)|,$\\[0.2cm]
 where $ \|\cdot \| $ and $|\cdot|$ denote matrix norm
and vector norm respectively.

It is easy to see that $~y(t)= (x^{1}(t)- x_{1}^{1}(t)) A_{1}
x(t)+ x_{1}^{1}(t) A_{1}(x(t)- x_{1}(t)). $ Then\\[0.2cm]
$|y(t)| \leq |(x^{1}(t)- x_{1}^{1}(t)) A_{1} x(t)| + |x_{1}^{1}(t)
A_{1}(x(t)- x_{1}(t))|. $ We have\\[0.2cm]
$|y(t)| \leq \|A_{1}\|( |x(t)|\cdot |x^{1}(t)- x_{1}^{1}(t))| +
|x_{1}^{1}(t)|\cdot|x(t)- x_{1}(t))|~ $ and using the
inequality\\[0.2cm]
$~|x^{1}(t)- x_{1}^{1}(t))|\leq |x(t)- x_{1}(t))|~$ one
obtains\\[0.2cm]
$(2)~~~|y(t)| \leq \|A_{1}\|( |x(t)|+ |x_{1}^{1}(t)|) |x(t)-
x_{1}(t))|.$\\[-0.2cm]

Similarly, we prove that\\[0.2cm]
$(3)~~~|z(t)| \leq \|A_{2}\|( |x(t)|+ |x_{1}^{2}(t)|) |x(t)-
x_{1}(t))|. $\\[-0.2cm]

According to $(2)$ and $(3)$, the relation $(1)$ becomes\\[0.2cm]
$|f(x(t)) - f(x_{1}(t))|\leq  (\|A\| + \|A_{1}\|( |x(t)|+
|x_{1}^{1}(t)|) + \|A_{2}\|( |x(t)|+ |x_{1}^{2}(t)|))
|x(t)-x_{1}(t)|.$\\[-0.2cm]

Replacing  $\|A\|= \|A_{1}\|=\|A_{2}\|= \sqrt{2}$, from the above
 we deduce that\\[0.2cm]
$(4)~~|f(x(t)) - f(x_{1}(t))|\leq  L |x(t)-x_{1}(t)|,~~~~~
\hbox{where}~~~L = \sqrt{2}(1+ 4 |x_{0}| +2 \delta )>0.$\\[-0.2cm]

 The inequality $(4)$ shows that $f(x(t))$
satisfies a Lipschitz condition. Using Theorems $1$ and $2$ in
\cite{chen}, it follows that $(3.4)$ has a unique solution.
\hfill$\Box$

 The equilibrium states of the fractional $5D$ Maxwell-Bloch
model $(3.3)$ are given as the union of the following two
families:\\[-0.3cm]
\[
E_{1}:=\{ e_{1}^{m,n}=(m, n, 0, 0, 0)\in {\bf R}^{5} |~m^{2} +
n^{2} \neq 0 \},~~~ E_{2}:=\{ e_{2}^{m}= (0, 0, 0, 0, m)\in {\bf
R}^{5} |~ m \in {\bf R}\}.
\]

\section{Stability study of fractional $5D$ Maxwell-Bloch
model}

We start with the study of stability of equilibrium states for the
fractional system $(3.3)$. Finally, we will discuss how to
stabilize the unstable equilibrium states of the system $(3.3)$
via fractional order derivative.

The Jacobian matrix of the system $(3.3)$ is\\[-0.2cm]
\[
J(x) = \left ( \begin{array}{ccccc}
0 & 0 & 1 & 0 & 0\\
0 & 0 & 0 & 1 & 0\\
x^{5} & 0 & 0 & 0 &  x^{1}\\
0 & x^{5} & 0 & 0 & x^{2}\\
-x^{3} & -x^{4} & - x^{1} & - x^{2} & 0\\
\end{array}\right ).\\[-0.1cm]
\]

\begin{Prop}
All equilibrium states of the fractional system  $(3.3)$ are
unstable.
\end{Prop}
{\it Proof.} {\bf Case $ e_{1}^{mn}\in E_{1}$.} The characteristic
polynomial of the matrix $J(e_{1}^{mn})$ is $~
p_{J(e_{1}^{mn})}(\lambda) = \det ( J(e_{1}^{mn}) - \lambda I)= -
\lambda^{3} (\lambda^{2} + m^{2} + n^{2}).$ Then the
characteristic roots of $ J(e_{1}^{mn}) $ are $
\lambda_{1}=\lambda_{2}=\lambda_{3} =0 $ and $ \lambda_{4,5}= \pm
i \sqrt{m^{2} + n^{2}}.$ Since the eigenvalues of $J(e_{1}^{mn})$
are at least one positive, by Proposition 2.1, it follows that $
e_{1}^{mn}$ is unstable.

{\bf Case $ e_{2}^{m}\in E_{2}$.} The characteristic polynomial of
the matrix $J(e_{2}^{m})$ is\\ $p_{J(e_{2}^{m})}(\lambda) = -
\lambda ( \lambda ^{2} - m)^{2}$ with characteristic roots
$~\lambda_{1}=0,\, \lambda_{2,3}= \sqrt{m}, \lambda_{4,5}= -
\sqrt{m} $ for $ m>0$ and $~\lambda_{1}=0,\, \lambda_{2,3}=
i\sqrt{-m}, \lambda_{4,5}= -i \sqrt{-m}$ for $ m<0$. Applying now
Proposition 2.1. it follows that $ e_{2}^{m}$ is unstable.

Similarly, it is easy to see that $ e_{0}$ is
unstable.\hfill$\Box$

The {\it controlled fractional $5D$ Maxwell-Bloch model}
associated to fractional $5D$ Maxwell-Bloch model $(3.3)$ is defined by:\\[-0.2cm]
\begin{equation}
\left\{ \begin{array}{ccl}
 D_{t}^{\alpha}{x}^{1}  &=&  x^{3}- k_{1}(x^{1}-x_{e}^{1})\\
  D_{t}^{\alpha}{x}^{2}  & =& x^{4}- k_{2}(x^{2}-x_{e}^{2})\\
 D_{t}^{\alpha}{x}^{3}  &=&  x^{1} x^{5}- k_{3}(x^{3}-x_{e}^{3}),~~~~~~~~~~~~~~~ \alpha\in (0,1),\\
  D_{t}^{\alpha}{x}^{4}  & =& x^{2} x^{5}- k_{4}(x^{4}-x_{e}^{4})\\
 D_{t}^{\alpha}{x}^{5}  &=& -( x^{1} x^{3}+ x^{2} x^{4})- k_{5}(x^{5}-x_{e}^{5})
\end{array}\right.\label{(4.1)}\\[-0.1cm]
\end{equation}
where $x_{e}$ represents an arbitrary equilibrium state of $(3.3)$
and $k_{i}\in {\bf R}, i=\overline{1,5}$ are non-negative
constants.

 The parameters $k_{i}, i=\overline{1,5}$ are feedback control
 gains which can make the eigenvalues of the linearized equation
 of the system $(4.1)$  satisfy one of the
 conditions of Proposition 2.1 or one of the fractional
 Routh-Hurwitz conditions \cite{ahme}, then the trajectories of
 the system $(4.1)$ asymptotically approaches the
 equilibrium state $x_{e}$.

The Jacobian matrix of the controlled fractional system $(4.1)$ is\\[-0.2cm]
\[
J(x,k) = \left ( \begin{array}{ccccc}
-k_{1} & 0 & 1 & 0 & 0\\
0 & -k_{2} & 0 & 1 & 0\\
x^{5} & 0 & -k_{3} & 0 &  x^{1}\\
0 & x^{5} & 0 & -k_{4} & x^{2}\\
-x^{3} & -x^{4} & - x^{1} & - x^{2} & -k_{5}\\
\end{array}\right ).\\[-0.1cm]
\]

\begin{Prop}
Let $ k_{i} >0$ for $ i=\overline{1,5}.$ Then the equilibrium
state $ e_{2}^{m}\in E_{2} $ of the controlled fractional system
$(4.1)$ is locally asymptotically stable for all $\alpha\in (0,
1],$  if one of the following conditions holds:

 $(1)~~|k_{1}-k_{3}|=|k_{2}-k_{4}|~$ and $~ m =
-\displaystyle\frac{1}{4}(k_{1}-k_{3})^{2}\neq 0;$

$(2)~~\max \{- \displaystyle\frac{1}{4}(k_{1}-k_{3})^{2},  -
\displaystyle\frac{1}{4}(k_{2}-k_{4})^{2} \} < m < \min \{ k_{1}
k_{3},  k_{2} k_{4}\};$

$(3)~~- \displaystyle\frac{1}{4}(k_{1}-k_{3})^{2} < m < \min \{-
\displaystyle\frac{1}{4}(k_{2}-k_{4})^{2}, k_{1} k_{3} \};$

$(4)~~-\displaystyle\frac{1}{4}(k_{2}-k_{4})^{2} < m  < \min \{
-\displaystyle\frac{1}{4}(k_{1}-k_{3})^{2}, k_{2} k_{4} \};$

$(4)~~  m <\min\{- \displaystyle\frac{1}{4}(k_{1}-k_{3})^{2}  , -
\displaystyle\frac{1}{4}(k_{2}-k_{4})^{2}\}.$
\end{Prop}
{\it Proof.} The Jacobian matrix of the system $(4.1)$ at the point $e_{2}^{m}$ is \\[-0.2cm]
\[
J(e_{2}^{m}, k)=\left ( \begin{array}{ccccc}
-k_{1} & 0 & 1 & 0 & 0\\
0 & -k_{2} & 0 & 1 & 0\\
m & 0 & -k_{3} & 0 &  0\\
0 & m & 0 & -k_{4} & 0\\
0 & 0 & 0 & 0 & -k_{5}\\
\end{array}\right ),\\[-0.1cm]
\]
whose characteristic polynomial $ p_{J(e_{2}^{m},k)}(\lambda) =
\det ( J(e_{2}^{m},k) - \lambda I)$ is\\[0.1cm]
 $p_{J(e_{2}^{m},k)}(\lambda) =- (\lambda +k_{5})[ \lambda ^{2} +
(k_{1} +k_{3})\lambda + k_{1} k_{3}- m][ \lambda ^{2} + (k_{2}
+k_{4})\lambda + k_{2} k_{4}- m]. $

Its characteristic roots are $~\lambda_{1}=
-k_{5},~\lambda_{2,3}=\displaystyle\frac{-(k_{1}+k_{3})\pm
\sqrt{(k_{1}-k_{3})^{2}+4m}}{2},$\\[0.1cm]
$\lambda_{4,5}=\displaystyle\frac{-(k_{2}+k_{4})\pm
\sqrt{(k_{2}-k_{4})^{2}+4m}}{2}.$ We denote:\\[0.1cm]
 $\Delta_{1} =
(k_{1}-k_{3})^{2} + 4m,~ \Delta_{2} = (k_{2}-k_{4})^{2} + 4m,~ u =
- \displaystyle\frac{1}{4}(k_{1}-k_{3})^{2},~ v =-
\displaystyle\frac{1}{4}(k_{2}-k_{4})^{2}.$

$(1)~$ Let $ \Delta_{1}= \Delta_{2}= 0.$ Then
$|k_{1}-k_{3}|=|k_{2}-k_{4}|$ and $ m =
-\displaystyle\frac{1}{4}(k_{1}-k_{3})^{2}\neq 0.$ The eigenvalues
$~\lambda_{1}=
-k_{5},~\lambda_{2,3}=\displaystyle\frac{-(k_{1}+k_{3})}{2},
\lambda_{4,5}=\displaystyle\frac{-(k_{2}+k_{4})}{2} $ are all
negative. Then $ |arg (\lambda_{i})| = \pi >
\displaystyle\frac{\pi}{2}\alpha ~$ for any $~\alpha\in (0,1] $
and so $ e_{2}^{m}$ is asymptotically stable.

$(2)$ Suppose that $ \Delta_{1} > 0 $ and  $ \Delta_{2}
> 0 .$ Then $ m >u $ and $m >v.$  We have $\lambda_{5} <0.$ The eigenvalues $\lambda_{i}, i=\overline{1,4}$
are all negative iff
 $u < m < k_{1} k_{3} $ and $ v < m < k_{2} k_{4}$. Hence, for
 $ \max\{u,v\} < m < \min\{ k_{1} k_{3}, k_{2} k_{4}\}$ it implies
that  $ e_{2}^{m} $ is asymptotically stable for $\alpha\in
(0,1].$

$(3)- (4)$ We suppose now that $ \Delta_{1} < 0 $ and  $
\Delta_{2}
>0 .$ It follows $ m < u $ and $ m  > v$. In this case the eigenvalues $\lambda_{4,5}$ are negative
iff  $ v < m < k_{2} k_{4}$. For $ \Delta_{1} < 0 $, we have
$\lambda_{2,3}=\displaystyle\frac{-(k_{1}+k_{3})\pm i
\sqrt{-\Delta_{1}}}{2}.$ Since $ Re(\lambda_{2,3})=
-\displaystyle\frac{1}{2}(k_{1}+k_{3}) < 0 $ we have
$|arg(\lambda_{2,3})|= \pi  > \displaystyle\frac{\alpha \pi}{2}$
for all $  0< \alpha \leq 1.$ Applying now Proposition 2.1 (i), we
can conclude that $e_{2}^{m}$ is asymptotically stable if $ v < m
< \min \{u, k_{2} k_{4}\} $ and $\alpha\in (0,1]$.

Similarly, we discuss the case $ \Delta_{1} > 0 $ and  $
\Delta_{2} <0 .$

$(5)$ Finally, we suppose $ \Delta_{1} < 0 $ and  $ \Delta_{2} < 0
.$ It follows $ m < u $ and $ m  < v$. We have
$\lambda_{2,3}=\displaystyle\frac{-(k_{1}+k_{3})\pm i
\sqrt{-\Delta_{1}}}{2},~\lambda_{4,5}=\displaystyle\frac{-(k_{2}+k_{4})\pm
i \sqrt{-\Delta_{2}}}{2}. $ Since  $ Re(\lambda_{i})< 0 $ for
$i=\overline{2,5},$ we have $|arg(\lambda_{i})|= \pi  >
\displaystyle\frac{\alpha \pi}{2}$ for all $  0< \alpha \leq 1.$
 By Proposition 2.1(i), $e_{2}^{m}$ is asymptotically stable iff  $ m <
\min \{u,v\} $ and $\alpha\in (0,1]$ . \hfill$\Box$

\begin{Ex}
{\rm By choosing the parameters $ k_{i}, i=\overline{1,5}$ that
satisfy one condition from Proposition 4.2, then the trajectories
of the controlled fractional  model are driven to the unstable
equilibrium point $e_{2}^{m}$. The parameters are selected as: $
k_{1}= k_{3}=\frac{1}{4},~ k_{2}=\frac{3}{2},~ k_{4}=\frac{2}{3},
~k_{5}>0.$ For $ m=-\frac{1}{8} $ we have $-\frac{25}{144} < m <
\min \{0,1\}.$ It follows that the stability condition $(4)$ of
Proposition 4.2 is achieved. This implies that, the trajectories
of the controlled fractional system $(4.1)$ converge to $ e_{2} =
(0,0,0,0, -\frac{1}{8})$ for any $\alpha\in (0,1].$ In this case
we have $\Delta_{1}=-\frac{1}{2} $ and $~\Delta_{2}=\frac{7}{36}$
for
  any $\alpha\in (0,1].$ The eigenvalues are
  $\lambda_{1}<0, \lambda_{2,3} = -\frac{1}{4} \pm
  i\frac{\sqrt{2}}{4}~$ and $~\lambda_{4,5} = -\frac{13}{12} \pm
  i\frac{\sqrt{7}}{12}<0.$}\hfill$\Box$
\end{Ex}

\begin{Prop}
 The equilibrium state $e_{0}$ of the controlled fractional system
$(4.1)$ is locally asymptotically unstable for $k_{i}
>0, i=\overline{1,5}$ and $\alpha \in (0,1)$.
\end{Prop}
{\it Proof.} The characteristic
polynomial of the Jacobian matrix $J(e_{0}, k)$ is\\
$p_{J(e_{0}, k)}(\lambda) = - \prod_{i=1}^{5}(\lambda + k_{i})$
with characteristic roots $~\lambda_{i}=-k_{i}$ for
$i=\overline{1,5}.$ Since\\[0.1cm]
 $ |arg(\lambda_{i}|=\pi >
\displaystyle\frac{\alpha \pi}{2}$ for $i=\overline{1,5}$, by
Proposition 2.1(i) it follows that $ e_{0}$ is locally
asymptotically stable.\hfill$\Box$

Let us we study the problem of stabilizing of the fractional
system $(3.3)$ at the equilibrium state $e_{1}^{mn}\in E_{1}$.

The Jacobian matrix of the system $(4.1)$ at the point $e_{1}^{m,n}$ is \\[-0.2cm]
\[
J(e_{1}^{mn}, k)=\left ( \begin{array}{ccccc}
-k_{1} & 0 & 1 & 0 & 0\\
0 & -k_{2} & 0 & 1 & 0\\
0 & 0 & -k_{3} & 0 &  m\\
0 & 0 & 0 & -k_{4} & n\\
0 & 0 & -m & -n & -k_{5}\\
\end{array}\right ).\\[-0.2cm]
\]

Its characteristic polynomial is $ p_{J(e_{1}^{mn},k)}(\lambda) =
- (\lambda + k_{1})(\lambda + k_{2}) P(\lambda)~
 $ with\\[-0.1cm]
\begin{equation}
P(\lambda) = \lambda^{3} + a_{1} \lambda^{2} + a_{2}\lambda +
a_{3},~~~~~\hbox{where}\label{(4.2)}
\end{equation}
\begin{equation}
\left \{\begin{array}{lcl}
a_{1} & = & k_{3} + k_{4} + k_{5}\\
a_{2} & = & k_{3} k_{4} + k_{3}k_{5} + k_{4} k_{5} + m^{2} + n^{2}.\\
a_{3} & = & k_{3} k_{4} k_{5} +  k_{3}n^{2} +  k_{4} m^{2}
\end{array}\right. \label{(4.3)}
\end{equation}

The eigenvalues of the characteristic equation are
$\lambda_{1}=-k_{1}, \lambda_{2} =-k_{2}$ and the roots
$\lambda_{3,4,5}$ of the equation $P(\lambda)=0.$

In this case we apply the fractional Routh-Hurwitz conditions
corresponding to polynomial $P(\lambda)$. The discriminant $D(P)$
of the polynomial $P(\lambda)$
is\\[-0.1cm]
\begin{equation}
D(P) = 18 a_{1}a_{2}a_{3} + a_{1}^{2} a_{2}^{2} - 4 a_{3}
a_{1}^{3} - 4 a_{2}^{3} - 27 a_{3}^{2}.\label{(4.4)}
\end{equation}

Because of the complexity of  $ D(P),$ we only consider the
following two situations:\\[0.2cm]
 $~~~~~(i)~~k_{i}>0 ~$ for
$i=\overline{1,5};~~~~~(ii)~~k_{1}>0,~ k_{2}
>0,~ k_{3}=k_{4} = b>0,~k_{5}= 0.$\\[-0.2cm]

In the above conditions we have $~a_{1}>0,~a_{2} >0,~a_{3}
>0,~a_{1}a_{2}-a_{3}>0.$

\begin{Prop}
Let $e_{1}^{mn}\in E_{1}$ the equilibrium state of the system
$(4.1)$.

 $(i)~$ Let $~ k_{i} >0$ for $i=\overline{1,5}.$\\[0.1cm]
$(1)~$ if $~D(P)>0$, then $e_{1}^{mn}$ is locally asymptotically
stable for $\alpha\in (0,1);$\\[0.1cm]
$(2)~$ if $~D(P)<0$, then  $e_{1}^{mn}$  is locally asymptotically
stable for $\alpha\in (0,\frac{2}{3}).$

 $(ii)~$ Let $~ k_{1}>0,~ k_{2}>0,~ k_{3}=k_{4} = b>0,~k_{5}=
 0.$\\[0.1cm]
$(1)~$ if $~b > 2 \sqrt{m^{2}+n^{2}},~$ then $e_{1}^{mn}$  is
locally asymptotically stable for $\alpha\in (0,1);$\\[0.1cm]
$(2)~$ if $~0 < b < 2 \sqrt{m^{2}+n^{2}},~$ then $e_{1}^{mn}$ is
locally asymptotically stable for $\alpha\in (0,\frac{2}{3}).$
\end{Prop}
{\it Proof.} $(i)(1)~$ From hypothesis we have $\lambda_{1} < 0 $
and $\lambda_{2} <0, a_{1} >0, a_{2} >0 $ and $a_{1}a_{2} >
a_{3}$.  When $D(P)>0$, the assertion $(i)$ of fractional
Routh-Hurwitz conditions (\cite{chen}, p. 704) is satisfied. But
Routh-Hurwitz conditions are the necessary and sufficient
conditions for the fulfillment of Proposition 2.1(i).  Then
$e_{1}^{mn}$ is asymptotically stable for any $\alpha\in (0,1)$.

$(2)~$ We have $\lambda_{1} < 0 $ and $\lambda_{2} <0, a_{1} >0,
a_{2} >0 $ and $ a_{3} > 0$. When $D(P)<0$, the assertion $(ii)$
of fractional Routh-Hurwitz conditions (\cite{chen}, p. 704) is
satisfied. As above, we deduce that $e_{1}^{mn}$ is asymptotically
stable for any $\alpha\in (0, \frac{2}{3})$.

$(ii)~$ For $k_{3}=k_{4}=b>0 , k_{5}=0$, we have $ a_{1}= 2b
>0,  a_{2}= b^{2}+ m^{2}+n^{2} >0 , a_{3}= b (m^{2} +n^{2}) $ and $ D(P)= (m^{2}+n^{2})^{2}[ b^{2}- 4(m^{2} +n^{2})].$
 Using the same manner as in demonstration of assertions $(i)$ we
 prove that $(ii)(1)$ and $(ii)(2)$ hold.\hfill$\Box$

\begin{Ex}
{\rm By choosing the parameters $ k_{i}, i=\overline{1,5}$ that
satisfy one condition from Propositions 4.4, then the trajectories
of the controlled fractional model are driven to the unstable
equilibrium point $ e_{1}^{mn}. $ If we select the parameters as
follows: $ k_{1}>0, k_{2}>0, k_{3}= k_{4}=0.5,~ k_{5}= 0 ~$ and
$~m^{2} + n^{2}= 0.25, $ then $a_{1} = 1,~
a_{2}=0.5,~a_{3}=0.125.$ Since
 $D(P)= -\frac{3}{64}<0 $ it follows that the stability condition $(2)$ of Proposition
 4.4 (ii) is achieved. This implies that, the trajectories of the system $(4.1)$ converge to
$ e_{1}^{mn} = (m,n,0,0,0)$ when $ m^{2}+n^{2}=0.25$ and
$\alpha\in (0,1)$. The eigenvalues are
  $~ \lambda_{3}= -0.5,~\lambda_{4,5}=-0.25 \pm 0.433 i.$
For example, substituting $ k_{1}= k_{2}= 1.2, k_{3}= k_{4}=0.5,~
k_{5}= 0 $ and $\alpha = 0.65 $ in $(4.1)$ we obtains that the
controlled fractional system is asymptotically stable at $ e_{1}=(
\displaystyle\frac{\sqrt{3}}{4},\displaystyle\frac{1}{4}, 0, 0,
0)$.}\hfill$\Box$\\[-0.2cm]
\end{Ex}

\section{Numerical integration of the fractional system $(4.1)$}

Consider the fractional differential equations\\[-0.2cm]
\begin{equation}
\left\{\begin{array}{lcl}
 D_{t}^{\alpha} x^{i}(t) & = &
F_{i}(x^{1}(t), x^{2}(t), x^{3}(t), x^{4}(t), x^{5}(t)),~~~t\in
(0, \tau),~\alpha\in (0,1)\\
x(0) &=& (x_{0}^{1}, x_{0}^{2}, x_{0}^{3}, x_{0}^{4}, x_{0}^{5})
\end{array}\right.\label{(5.1)}
\end{equation}
where $~F_{1}(t) = x^{3}(t)- k_{1}(x^{1}(t) - x_{e}^{1}),~~
F_{2}(t) =  x^{4}(t)- k_{2}(x^{2}(t) - x_{e}^{2}),$\\[0.1cm]
$F_{3}(t) = x^{1}(t) x^{5}(t)- k_{3}(x^{3}(t) - x_{e}^{3}),~~
F_{4}(t)= x^{2}(t) x^{5}(t)- k_{4}(x^{4}(t) -
x_{e}^{4}),$\\[0.1cm]
$F_{5}(t) = -( x^{1}(t) x^{3}(t) + x^{2}(t) x^{4}(t) )-
k_{5}(x^{5}(t) - x_{e}^{5}).$

 Since the function $ F(t)= ( F_{1}(t), F_{2}(t), F_{3}(t), F_{4}(t), F_{5}(t))
 $is continuous, the initial value problem $(5.1)$ is equivalent to the nonlinear Volterra integral equation \cite{difo},
 which is given as follows:\\[-0.2cm]
\begin{equation}
x^{i}(t)~=~ x_{0}^{i}  +
\displaystyle\frac{1}{\Gamma(\alpha)}\int\limits_{0}^{t}(t-s)^{\alpha
- 1} F_{i}( x^{1}(s), x^{2}(s), x^{3}(s), x^{4}(s), x^{5}(s)) d s,
~~~~~i=\overline{1,5}.\label{(5.2)}
\end{equation}

Diethelm et al. have given a predictor-corrector scheme
\cite{difo}, based on the Adams-Bashforth-Moulton algorithm to
integrate the equation $(5.2)$. We apply this scheme to the
controlled fractional system $(5.1)$. For this, let $ h =
\displaystyle\frac{\tau}{N}, ~ t_{n}=n h~$ for
$ n = 0,1,\ldots, N.$ We use the following notations:\\[0.1cm]
$ x^{i}[n] = x^{i}(n h),~~ x_{p}^{i}[n] = x_{p}^{i}(n h),~~
F_{i}[n] = F_{i}(x[n]), ~~ F_{i,p}[n]= F_{i}(x_{p}[n])~$ for $~i=\overline{1,5}.$\\

 The controlled fractional system
$(5.1)$ can be discretized as
follows:\\[-0.2cm]
\begin{equation}
\left\{\begin{array}{lcl}
x_{p}^{i}[n+1]
 &= &  \displaystyle\frac{h^{\alpha}}{\alpha\Gamma(\alpha
)} \sum\limits_{j=0}^{n}b[j,n+1]
F_{i}[j] \\[0.2cm]
x^{i}[n+1] & = & x_{0}^{i} +
\displaystyle\frac{h^{\alpha}}{\Gamma(\alpha +
2)}(\sum\limits_{j=0}^{n}a[j,n+1] F_{i}[j]+  F_{i,p}[n+1]
),\label{(5.3)}
\end{array}\right.
\end{equation}
where $~ i=\overline{1,5}~$ and:\\[-0.3cm]
\begin{equation}
\left\{ \begin{array}{lcl}
a[0,n+1] & = & n^{\alpha + 1}-( n - \alpha)(n + 1)^{\alpha}\\
a[j,n+1] & = & (n - j + 2)^{\alpha + 1}+( n - j)^{\alpha + 1}- 2
(n - j + 1)^{\alpha + 1},~~~j=\overline{1,n}\\
b[j,n+1]& =& (n+1-j)^\alpha - (n
-j)^{\alpha},~~~j=\overline{0,n}.\label{(5.4)}
\end{array}\right.
\end{equation}

The above scheme given by the relations $(5.3)$ and $ (5.4)$  is
called the {\it Moulton- Adams algorithm for controlled fractional
system $(5.1)$} (see for details \cite{difo}).

The error estimate for the algorithm described by $(5.3)$ and $
(5.4)$ is\\[-0.2cm]
\[
\max_{ 0\leq j \leq N}\{~ x^{i}[j] - x_{p}^{i}[j]) |
i=\overline{1,5}~\}= O(h^{\alpha + 1}).\\[-0.2cm]
\]

Applying the algorithm $(5.3)-(5.4)$, the fractional system
$(5.1)$ is numerically integrated for $\alpha = 0.65,~ k_{1}=
k_{2}= 1.2, k_{3}= k_{4}=0.5, k_{5}= 0 $ and $ x_{e}=(
\displaystyle\frac{\sqrt{3}}{4},\displaystyle\frac{1}{4}, 0, 0,
0)$ (see Example 4.2). For this, we consider $~h = 0.01,
\varepsilon= 0.01,  N = 500, t = 502 $ and the initial conditions
$~x^{1}(0)=\varepsilon + \displaystyle\frac{\sqrt{3}}{4},
  x^{2}(0)=\varepsilon + \displaystyle\frac{1}{4}, x^{3}(0)=x^{4}(0)=x^{5}(0)=\varepsilon.$

Using the software Maple $11$, the orbits $(n, x^{i}(n)),
i=\overline{1,5}~$ of system $(5.1)$ are represented in the
figures Fig. 1-5.

\begin{center}
\begin{tabular}{ccc}
\includegraphics[width=5cm]{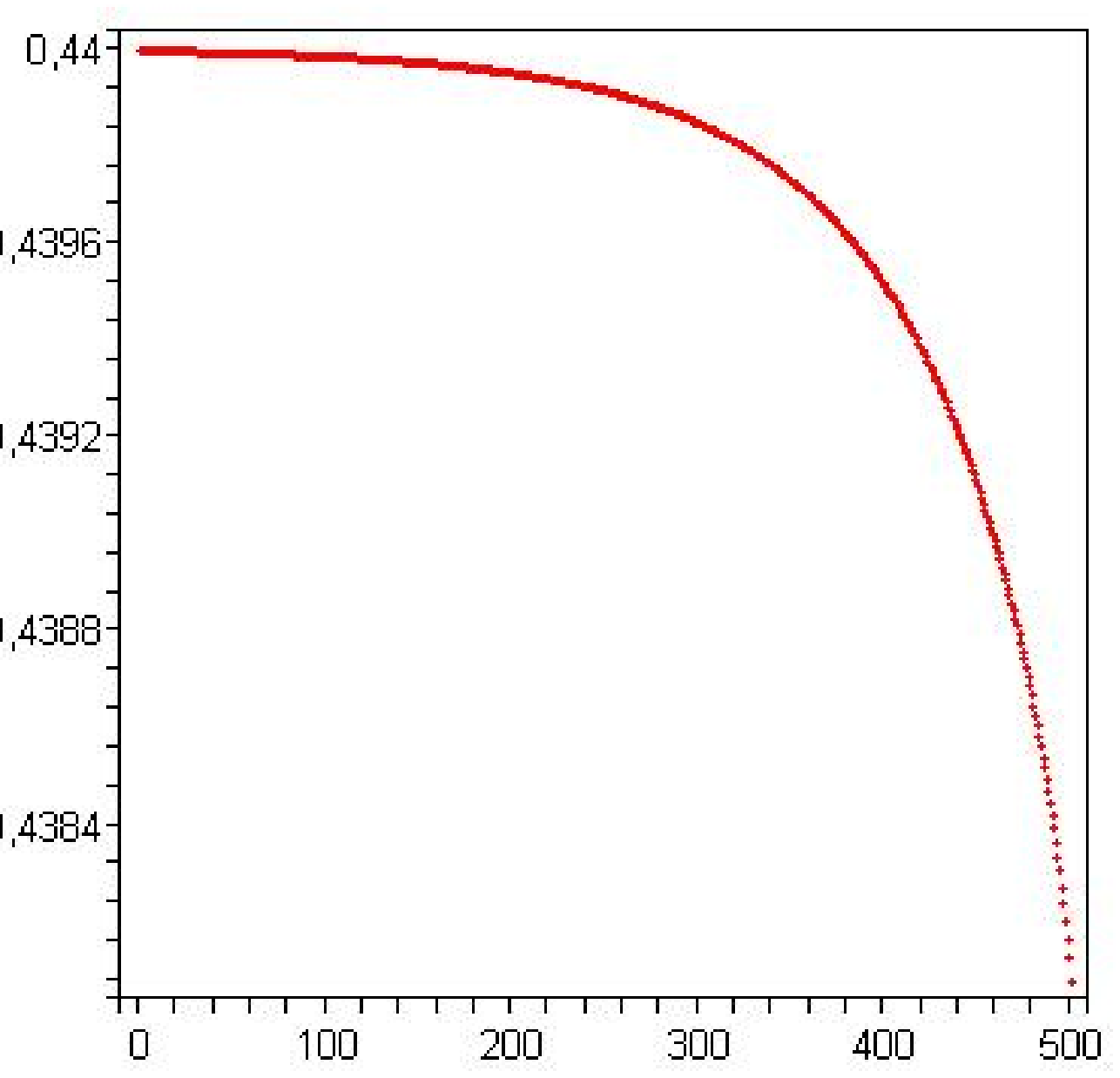}&
\includegraphics[width=5cm]{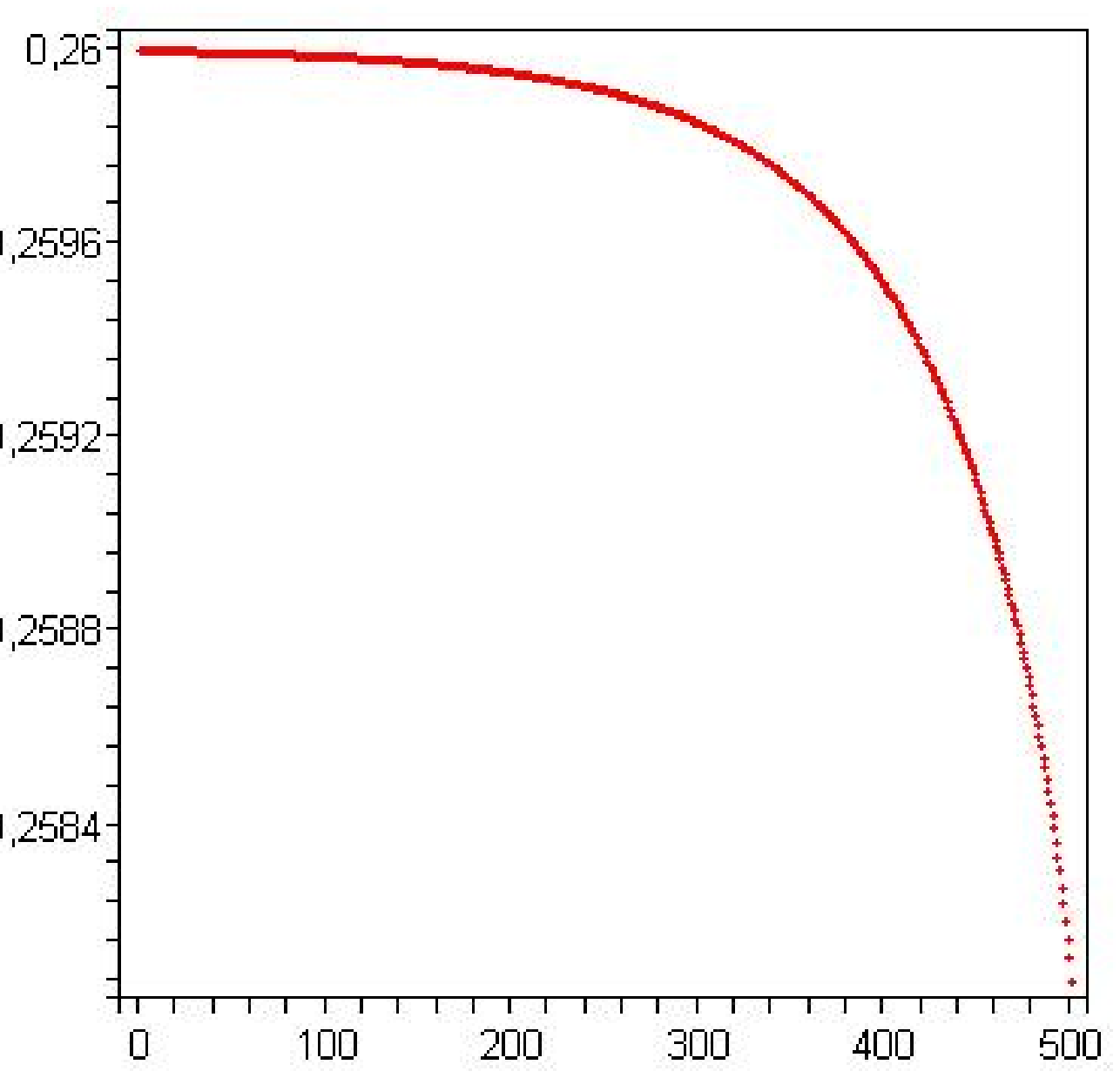}&
\includegraphics[width=5cm]{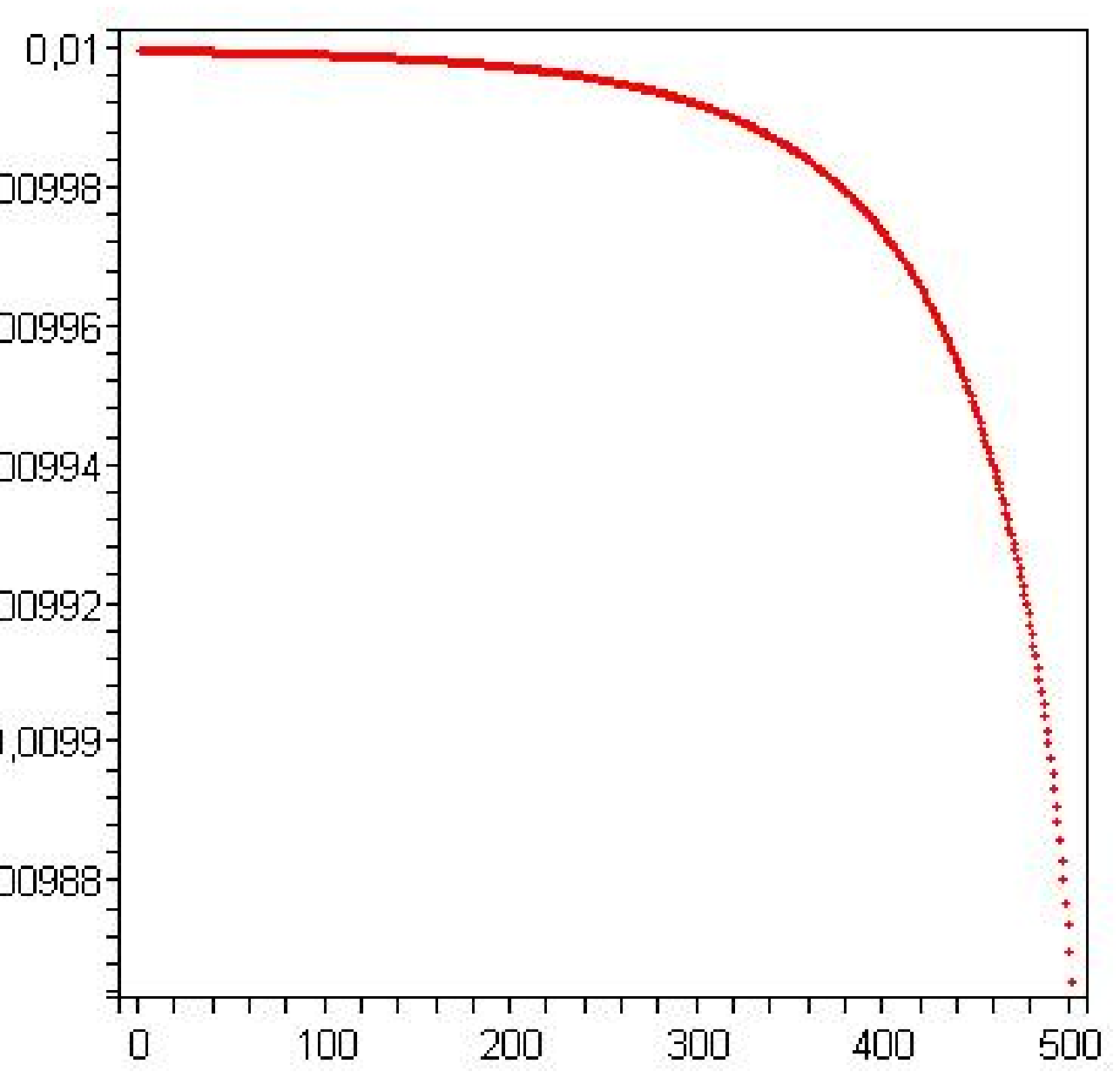}
\end{tabular}
\end{center}

\begin{center}
\begin{tabular}{ccc}
$~~~~~$ {\bf Fig.1.} $~(n, x^{1}(n))~~~~~~~~~~~~~~~~~$ & {\bf
Fig.2.} $~(n, x^{2}(n))~~~~~~~~~~~~~~~~~$ & {\bf Fig.3.} $~(n,
x^{3}(n)$\\
\end{tabular}
\end{center}

\begin{center}
\begin{tabular}{ccc}
\includegraphics[width=5cm]{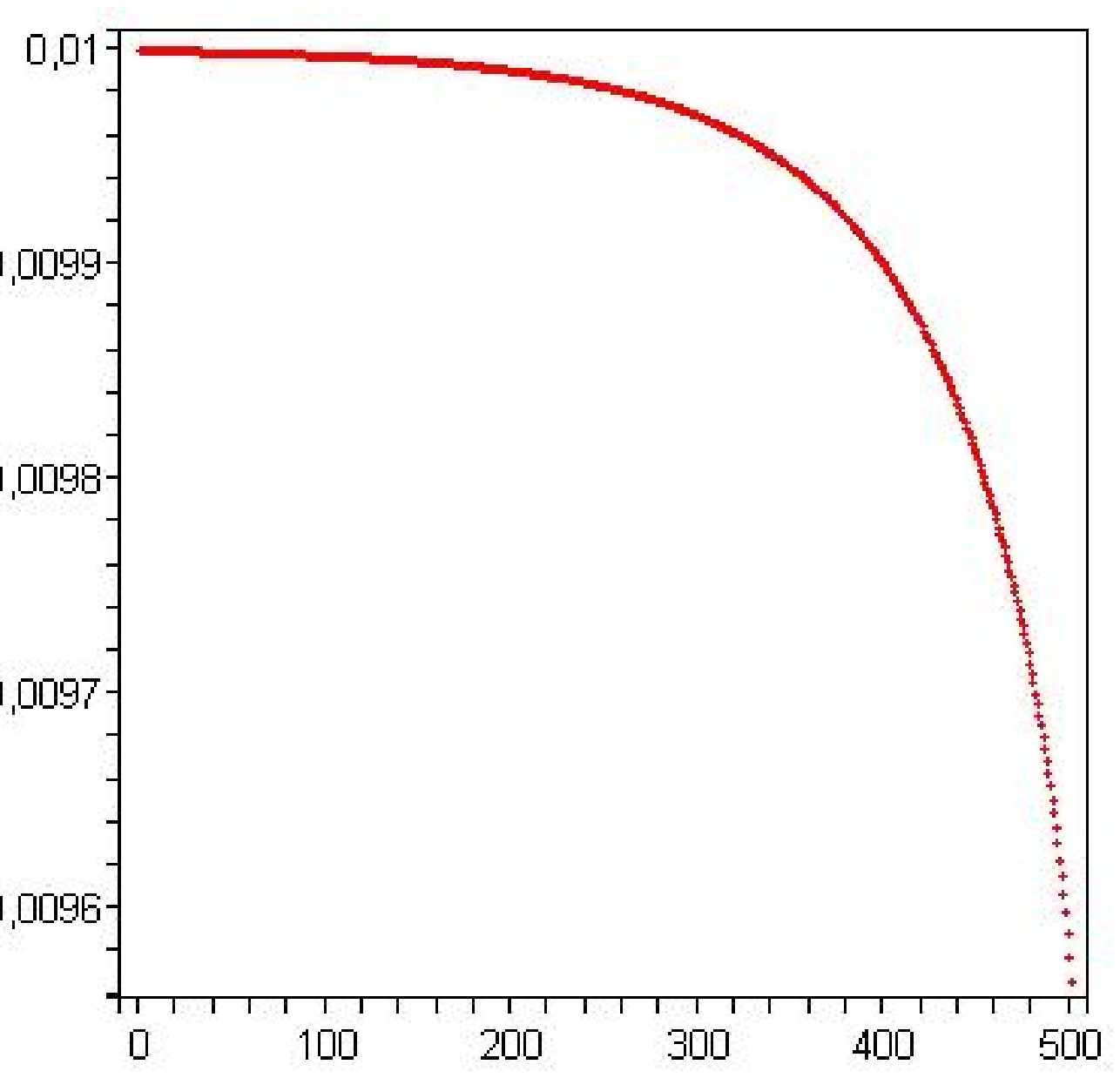}&
\includegraphics[width=5cm]{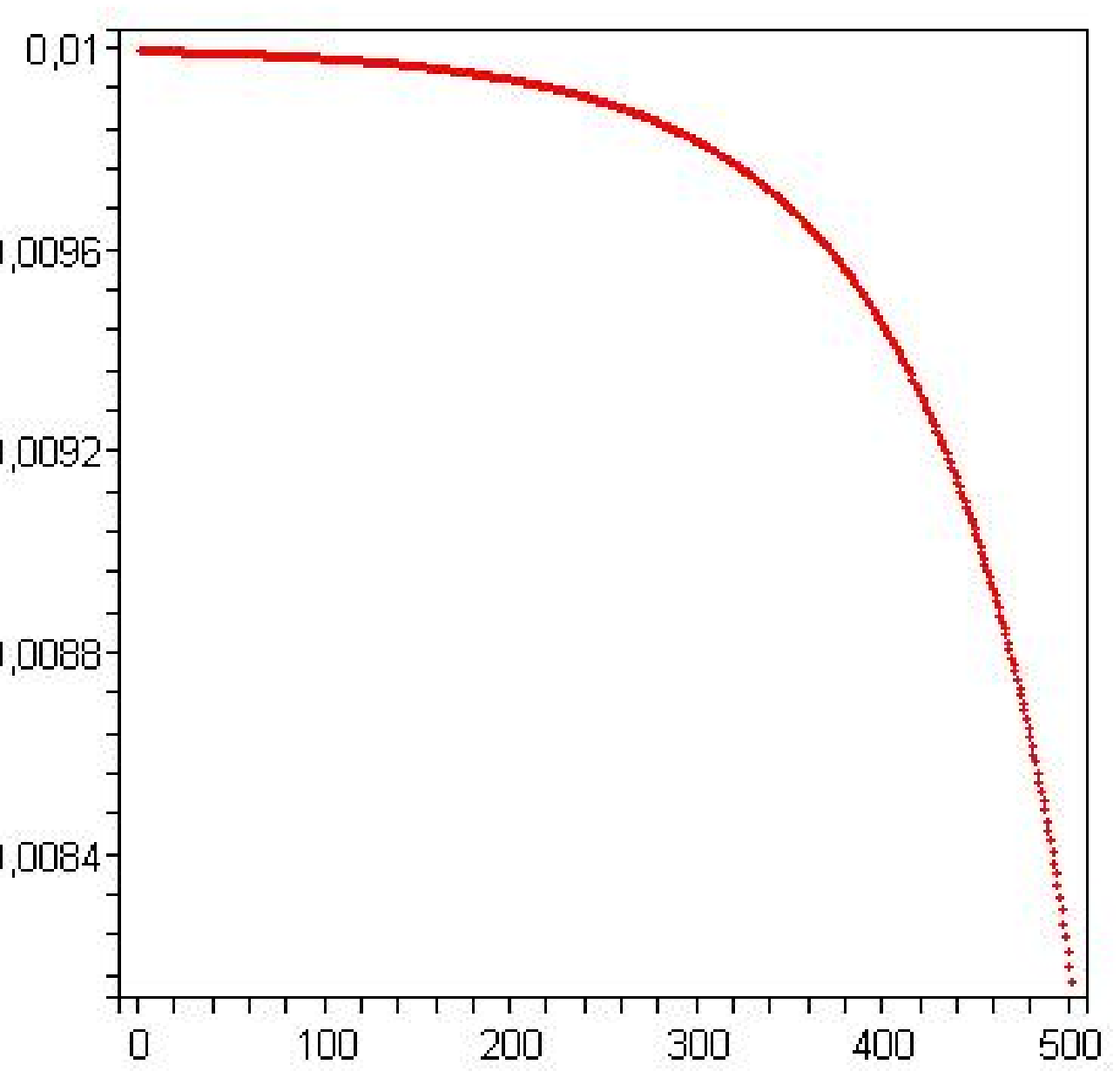}&
\end{tabular}
\end{center}

\begin{center}
\begin{tabular}{ccc}
$~~~~~~~~~~~~$ {\bf Fig.4.} $~(n, x^{4}(n))~~~~~~~~~~~~~~~~$ & {\bf
Fig.5.}
$~(n, x^{5}(n))~~~~~~~~~~~~~~~~~$ & \\
\end{tabular}
\end{center}

The numerical simulations show the validity of the theoretical
analysis.

{\bf Conclusions.} The dynamics of the fractional $5D$ Maxwell-Bloch
model $(3.3)$ was investigated in this paper. The analysis of the
stability of equilibrium states for the controlled fractional $5D$
Maxwell-Bloch model $(4.1)$ was studied. Finally, the numerical
integration and numerical simulation for the fractional system
$(4.1)$ are given.

 {\bf Acknowledgments.} The author has very grateful to be
reviewers for their comments and suggestions.\\[-0.3cm]

Author's adress\\

West University of Timi\c soara,\\
Seminarul de Geometrie \c si Topologie,\\
 Department of Educational Sciences,\\
4, B-dul V. P{\^a}rvan, 300223, Timi\c{s}oara, Romania.\\
E-mail: ivan@math.uvt.ro\\

 \end{document}